\magnification 1200

\def\intro{{{1.1. }}}
\def\plane{{{1.2. }}}

\def\GITSRP{{{2.1. }}}
\def\GIT{{{2.2. }}}
\def\extend{{{2.3. }}}
\def\extendlemma{{{2.3.1 }}}
\def\mumford{{{2.4. }}}

\def\CLASSIC{{{3.1. }}}
\def\table{{{3.2. }}}
\def\irreducible{{{3.3. }}}
\def\irreduciblelemma{{{3.3.1 }}}
\def\reducible{{{3.4. }}}
\def\tacnodal{{{3.5. }}}
\def\key{{{3.5.1 }}}

\def\split{{{4.1. }}}
\def\splitpr{{{4.1.1 }}}
\def\trinodal{{{4.2. }}}
\def\trinodalpr{{{4.2.1 }}}

\def\binodal{{{5.1. }}}
\def\binodalpr{{{5.1.1 }}}
\def\nodalpr{{{5.1.2 }}}
\def\smooth{{{5.2. }}}
\def\main{{{5.2.1 }}}

\font\tenmsb=msbm10
\font\sevenmsb=msbm10 at 7pt
\font\fivemsb=msbm10 at 5pt
\newfam\msbfam
\textfont\msbfam=\tenmsb
\scriptfont\msbfam=\sevenmsb
\scriptscriptfont\msbfam=\fivemsb
\def\Bbb#1{{\fam\msbfam\relax#1}}

\def\sqr#1#2{{\vcenter{\vbox{\hrule height.#2pt \hbox{\vrule width.#2pt
height #1pt \kern #1pt \vrule width.#2pt}\hrule height.#2pt}}}}
\def\qed{$\  \sqr74$}

\def\mapdown#1{\Big\downarrow \rlap{$\vcenter{\hbox{$\scriptstyle#1$}}$}}

 \def\la{\longrightarrow}

 \def\ni{\noindent}
 \def\cl{\centerline}
\def\df{\noindent {\bf Definition.\  }}
 \def\rk{\noindent {\it Remark.\  }}
 
 \def\pf{\noindent {\it Proof.\  }}

\def\d{\delta }

\def\g{\gamma}
\def\t{\theta}
\def\T{\Theta}

\font\gothicf=eufm10
\font\sgothic=eufm7
\font\ssgothic=eufm5
\textfont5=\gothicf
\scriptfont5=\sgothic
\scriptscriptfont5=\ssgothic

\def\C{{\Bbb C}}

\def\P{{\Bbb P}}
\def\Pl{\P ^2}
\def\Pld{\P ^{2*}}

\def\X{{\cal X}}
\def\Y{{\cal Y}}
\def\O{{\cal O}}
\def\Z{{\cal Z}}
\def\L{{\cal L}}

\def\Vs{V^0}

\cl{\bf{RECOVERING PLANE CURVES  FROM THEIR BITANGENTS}}

\

\cl{Lucia Caporaso and Edoardo Sernesi}

\

\

\cl{1.INTRODUCTION}

\

\ni
{\bf{\intro }}
It is well known that a general  plane curve of degree $d$ has ${1\over 
2}d(d-2)(d^2-9)$ distinct bitangent lines. The first (and most) interesting case 
is that of a smooth plane quartic $X$, whose configuration of $28$ 
bitangents we shall denote by $\t(X)$, to highlight the  
correspondence with the odd theta-characteristics of $X$.

The properties of  $\t(X)$ have been extensively studied by the classical
geometers since the time of Riemann  
(Aronhold,  Cayley, Hesse, Pl\"ucker, Schottky, Steiner, Weber just to quote a few). 
Such an exceptional interest is not only due to its  rich geometry,
 but also and especially to the connection of this configuration with the
classical theory of theta functions.

Some of the features  of $\t(X)$ discovered classically have to do with the 
possibility of recovering the curve $X$ from various data related to $\t(X)$. The
most celebrated of these results is  due to Aronhold: he discovered  for
every nonsingular quartic $X$ the existence of $288$ $7$-tuples of bitangents,
called {\it Aronhold systems}, characterized by the condition that for any three
lines of such a system the six  points of contact are not  on a conic; he showed that
the full  configuration $\t(X)$ and the curve $X$ can be reconstructed  from any one of
these Aronhold systems. For an account of his construction we refer  to  [E-Ch],
p. 319, or [K-W], p. 783.

Notice that the definition of Aronhold system requires knowing not only $\t(X)$, but
also the contact points between $X$ and its bitangents.
And since there is no known way
to recognize the Aronhold systems only knowing  $\t(X)$,    
the result of Aronhold leaves still unanswered the question below, to which this paper is
devoted:

\

\ni
{\bf Question} (1st form): {\it Can the curve $X$ be recovered from  
its bitangent lines?}

\

Or: if $X$ and $Y$ are two nonsingular plane quartics such that 
$\t(X) = \t(Y)$, does it follow that $X=Y$?

Interest in this question is also related to the
classical Torelli theorem. Consider  the jacobian
variety $J(X)$ of $X$, and let  $\T _X\subset J(X)$ be a (symmetric) theta divisor. Then
the bitangent  lines are the images, under the Gauss map
$\T _X\to \Pl$,  of the   points of $\T _X$ which are $2$-torsion points of 
$J(X)$ (see [ACGH] for details); 
a positive answer to the above question would therefore imply the possibility
of recovering $X$ from  the Gauss images of finitely many points of $\T _X$, thus 
giving a refined and sharper version of the classical Torelli theorem.

Recalling that a nonsingular plane quartic is a canonical curve, our question can be therefore more
fruitfully restated as follows:

\ 

\ni
{\bf Question} (2nd form) : {\it Can a nonsingular, non hyperelliptic curve $X$ of genus 3 be
reconstructed from the Gauss images of the (smooth) $2$-torsion points of the
theta-divisor $\T _X\subset J(X)$?}

\ 

The question, in the second form, can be asked for curves of any genus 
$g\ge 3$. This was in fact the main motivation for this work: we
plan to come back to such a generalization in a subsequent paper. 

\

\ni
{\bf \plane }
On the other hand the question, in the first form, can be asked for plane  curves  of
any degree $d \ge 4$. It is a remarkable fact that for curves of degree 
$d \ge 5$ the answer is affirmative for simple reasons; in other words, every
general plane curve of degree $d \ge 5$ is uniquely determined by its 
bitangent lines. In fact 
 each of the ${1\over 2}d(d-2)(d^2-9)$ bitangents  corresponds to a node of the dual curve 
(which has degree
$d(d-1)$). Thus, if two nonsingular curves $Y_1$ and $Y_2$ of degree 
$d$ have
the same bitangents,
 their dual curves
have all the corresponding double points  
in common, and these contribute with $4[{1\over 2}d(d-2)(d^2-9)]$ to the degree of
the intersection of the dual curves. This contradicts B\'ezout if $d \ge 5$. 

This argument fails for $d=4$, and no other argument of an elementary
nature seems to work.
Our approach has therefore been different. It is not too difficult to show (see Section 4)
that
certain singular quartics are uniquely determined by their generalized bitangent
lines, i.e. they have the {\it theta-property} (see \extend for the definitions). 
Then, combining this with a 
degeneration argument, we have been able to give a positive answer to our Questions
at least in a weak form. Our main theorem is the following (\main ).

\

\proclaim THEOREM. The general plane quartic  is uniquely
determined by  its 28 bitangent lines.

We can actually prove that the same is true for a  singular curve, 
 so that a general quartic having $\d$ nodes
($\d=0,\ldots, 4$) has the theta-property (see sections 4 and 5 for details).

We use a degeneration technique; after having 
generalized in a natural way the definition of $\t (X)$ to certain singular quartics (see sections 2 and 3), we show that
every irreducible quartic $X$ having
exactly $3$ nodes is uniquely determined by its own $\t (X)$ (\trinodalpr ). 
Then we consider suitable families of smooth 
curves specializing to trinodal curves. To obtain our result we need to eliminate ``bad" limits
(such as non-reduced curves). To do that we apply  Geometric Invarian Theory   to the
action of $PGL(3)$  on two different spaces: the space of all plane quartics and the space $Sym ^{28}(\Pld )$, where the
configurations of bitangents live.

\

\cl
{2. PRELIMINARIES.} 

\

\ni
{ \bf \GITSRP}
We work over $\C$.

Let us here recall some fundamental facts of Geometric Invariant Theory, including the so-called
GIT-semistable replacement property.
Given  a projective scheme $H$ (smooth, irreducible for simplicity)  over which a reductive
group
$G$  (for our purposes, $G=PGL(3) $ or  $G=SL(3) $) acts in a linear way,
let
$H^{SS}$ be the open subset of GIT-semistable points. Then there exists 
a  quotient (in a suitable sense) morphism
$$
q:H^{SS}\la H^{SS}//G
$$
with $H^{SS}//G$ a projective scheme.
The complement of $H^{SS}$ in $H$ is the locus of  GIT-unstable points.
If $x$ and $y$ are GIT-semistable points in $H$, then $q(x)=q(y)$ if and only if 
$$
\overline{O_G(x)}\cap \overline{O_G(y)}\cap H^{SS} \neq \emptyset
$$
where $
 O_G(x)$ denotes the orbit of $x$ via $G$.
The GIT-semistable replacement property is just an application of the existence of such a projective
quotient.
 Let $T$ be a smooth curve and $t_0\in T$ be a point. 
Let $\phi : T - \{t_0\} \la H^{SS}$ be a regular morphism. Then there exists a finite covering $\rho:T'\la T$,
ramified only over $t_0$ (let $t'_0=\rho^{-1}(t_0)$) such that
 there exist morphisms $\psi :T' \la H^{SS}$ and $\g :T'-\{t'_0\}\la G$ with the property that
$\forall t'\neq t'_0$
we have $\phi (\rho(t'))^{\g (t') }= \psi (t')$.
Moreover,  $\psi (t'_0)$ can be chosen arbitrarily in a (uniquely determined) fiber of $q$, in particular,
one can assume that it
has finite stabilizer.

\

\ni
{\bf \GIT }
The space of all plane quartics will be identified with $\P ^{14}$. The group $G=PGL(3)$
of  automorphisms of $\Pl$ acts on
it in a natural way. Such an action is linear in the sense of GIT; we shall denote by $V'\subset \P
^{14}$ the open subset of quartics that are semistable with respect to the action of $G$.
It is well known ([GIT]) that $V'$ is made of   quartics having at most double points as
singularities (not all of them, see below).
Thus the only non-reduced quartics that are GIT-semistable are double conics. 
The (closure of the) set of all double conics in $\P ^{14}$ is a closed subset of dimension 5;
we let $V$  be the complement in $V'$ of such closed subset.  Recall also that $X\in V$ is
GIT-stable if and only if $X$ is reduced and has no tacnodes. In particular, if $X$ has only
ordinary nodes and ordinary cusps it is GIT-stable (we refer to [GIT], chapter 4,
section 2).
Notice that $V$ is open and  contains  the set $\Vs$ of all smooth plane quartics.

\

\ni
{\bf \extend}
Let $X$ be a plane quartic. We say that a line $L$ is a bitangent of $X$ if $L$ intersects $X$
in smooth points and if the scheme $X
\cap L$ is  everywhere nonreduced.

 Recall  that every smooth plane quartic has exactly 28 distinct bitangent
lines.  Denote $P_n:=Sym^{n}(\Pld )$; we can  define a map
$$
\t :\Vs \la P_{28}
$$
such that for every $X\in \Vs$, $\t (X)$ is the set of 28 bitangents of $X$.
We shall  use the notation $\t (X)$ also to indicate the  corresponding plane curve of degree
$28$. We call
$\t (X)$ the {\it theta-curve} of
$X$.

Our goal now
is to extend such a definition to all curves in $V$.
First we need the following

\

\df
Let $X$ be a plane quartic. A line $L$ is called a {\it theta-line} of $X$ if either $L\subset
X$, or if
$X\cap L$ is everywhere non-reduced. Let $X$ be reduced and  not containing $L$;  if
$i $ is a non-negative
integer and 
$L$ contains exactly
$i$ singular points of
$X$, we shall say that $L$ is a theta-line of type $i$.

\

Some examples: the bitangent lines of a smooth curve are theta-lines of type zero.
If $X$ has one node (one cusp) and no other singularities, then (by  the Hurwitz formula applied to the projection from the
singular point) there exist exactly $6$ distinct lines passing through the node (cusp) and tangent to $X$. These
$6$ lines are theta-lines of type $1$. 

We shall see more examples later.
Now we prove the following simple:

\proclaim Lemma \extendlemma. There exists a natural extension of the morphism $\t$ to the whole of $V$, such
that for every $X\in V$, all components of  $\t (X)$  are theta-lines of $X$.

\rk In particular,  if $X_0$ is any singular curve in $V$ and $\X \la T$ is any
deformation
of $X_0$ to smooth curves, then the limit $\T _0$ of the theta curves $\t(X_t)$ of the smooth fibers does not
depend on the choice of the deformation.

\

\pf
Let $X\in V$; then (see \GIT ) $X$ is reduced and has no point of multiplicity greater
than 2.
Thus, it is easy to see that the set of all theta-lines of $X$ is finite.

Now consider the incidence correspondence
$$
J^0\subset \Vs \times P_{28}
$$
defined by $J^0 =\{ (X; \t (X)):\  X\in \Vs \}$, or
$$
J^0=\{ (X;L_1,....,L_{28}):\  L_i\neq L_j,\  L_i\cap X\  
\rm{is \  everywhere \  nonreduced\  } \forall i=1,....,28\}.
$$
Let $J$ be the closure of $J^0$ in $V\times P_{28}$ and let $\rho:J\la V$ be the projection.
If $X$ is any curve in $V$ and $\T \in \rho ^{-1} (X)$, then $\T$ is a set of 28 lines
(not necessarily distinct) and such that if $L$ is a line appearing in $ \T$, $L$ is 
a theta-line of $X$.
Since there are only finitely many theta-lines for any given $X\in V$ and since $\deg \T =28$ we
conclude that $\rho$ has finite fibers.
Of course, $\rho$ is one-to-one on the open subset $V^0$ of smooth curves; moreover
 $J$ is irreducible (because $\Vs $, and hence $J^0$, is irreducible). Therefore, by the Zariski
connectedness principle ($V$ is smooth) we conclude that $\rho $ is one-to-one everywhere. Finally,  for
every $X\in V$ we can define $\t (X)$ by the rule
$\rho ^{-1}(X)=\{ (X,\t (X))\}$ and we are done.
\qed

\

\df
For any curve $X$ in $V$, the curve represented by $\t (X)$ (defined in the above Lemma) will be called
the {\it theta-curve} of $X$.

\

Section 3 is devoted to a precise description of the theta-curve of  a singular quartic.

\

\df
Let $X$ be a curve in $V$, and let $\T =\t (X)$ be its theta-curve. If $X$ is the only curve in
$V$ having $\T$ as theta-curve (i.e. if $Y\in V$ is such that $\t (Y)=\T$, then $Y=X$),
then we shall say that $X$ has the {\it theta-property}.

\

\ni
{\bf \mumford}
Our goal is to show that a general quartic has the theta-property. To do that we shall consider the
action of $G$ on  $P_n:=Sym^{n}(\Pld )$. David Mumford, in [GIT] 4.4, gives necessary and
sufficient conditions for a point in
$P_n$ to be GIT-stable, semistable and unstable with respect to the natural action of  $G$. His criterion is
the following: $\Sigma \in P_n$ is semistable if and only if the plane curve of degree $n$ corresponding to
$\Sigma$ contains no point of multiplicity greater than $2n/3$ and no line of multiplicity greater than
$n/3$. In particular a point $\Sigma$ in $ P_{28}$ is
GIT-unstable if and only if the corresponding plane curve of degree $28$  
has either a  point of multiplicity at least
$19$ or a line of multiplicity at least $10$. Otherwise, if $\Sigma$ has at most points
of multiplicity $18$ and lines of multiplicity $9$, $\Sigma$ is GIT-stable (there are
no strictly semistable points in this situation).

\

\

\cl{3. THE STRUCTURE OF THE THETA CURVE OF A SINGULAR QUARTIC}

\

\ni
{ \bf\CLASSIC}
Let $X$ be an irreducible quartic in $V$. Then $\t(X)$ consists of 28
lines of which $b_i$ are of type $i$, $i=0,1,2$,
and each of them  appears with a  multiplicity  to be computed.

Recall that for an irreducible plane curve $C$ having only nodes,
 ordinary cusps and ordinary tacnodes as singularities, and only 
 ordinary flexes, the degree 
$d$, the class (i.e. the degree of the dual curve) $m$, and the numbers 
$\d $, $\kappa$, $\tau$ $f$, $b$ of nodes, cusps, tacnodes, flexes and bitangents  (respectively)
are related by the classical Pl\"ucker formulas: 
$$
m = d(d-1) - 2\d  - 3\kappa - 4\tau
$$
$$
d = m(m-1) - 2b - 3f -  4\tau
$$
$$
f = 3d(d-2) - 6\d  - 8\kappa - 12 \tau 
$$
(see  [W]). From these formulas one 
easily computes the following expression for $b$ in terms of $d$, $\d $, 
$\kappa$ and $\tau$:
$$
b = N_d - (d+2)(d-3)(2\d  + 3 \kappa + 4 \tau) + 2\d (\d  - 1+ 4\tau) 
+6\kappa(\d+2\tau) +
{9\kappa(\kappa-1)\over 2}+ 2\tau(4\tau-3)
$$
where
$$
 N_d = {1\over 2} d(d-2)(d^2-9)
 $$
 is the number of bitangents of a nonsingular curve of degree $d$. 
 
 Let's specialize to the case of irreducible quartics.
 The above
formula computes 
 the number $b_0$ of theta-lines of type $0$, giving:
$$
b_0 = 28 - 6(2\d+3\kappa+ 4 \tau)  + 2\d (\d  - 1+ 4\tau) 
+6\kappa(\d+2\tau) +
{9\kappa(\kappa-1)\over 2}+ 2\tau(4\tau-3)
$$
Obviously
$$
b_2 = {\d+\kappa+ \tau  \choose 2}
$$
The number $b_1$ can be computed as follows. 
From the Hurwitz formula it follows that through each double point 
the number of theta-lines is
$ 2(3-\d-\kappa-2\tau)+2$ (just consider the degree-$2$ map to $\P ^1$ given by projecting from the double point). Among
these one has the
 theta-lines  joining the double point with the cusps different from it,
if any, which are of type $2$, and therefore already counted. 
This simple rule allows us to compute $b_1$ case by case.  

Notice that these computations are valid even if we allow the 
irreducible quartic $X$ to have hyperflexes (i.e. nonsingular points 
where the tangent line has multiplicity of intersection equal to four).
In fact  each such hyperflex diminishes by two the number of ordinary 
flexes. 
At the same time it corresponds to a triple 
point of the dual curve which diminishes by eight its class (see [S], \S 21, for details).
All this 
accounts for an obvious modification in the second of the above 
Pl\"ucker formulas. Since the first and third formula obviously remain 
unchanged,  the final expression of $b_{0}$ stays the same.

\

\ni
{ \bf \table}
The following table lists the resulting numbers
 for all possible cases of irreducible quartics in $V$:
$$
\matrix{
(\d,\kappa,\tau) & b_0 & b_1 & b_2 \cr
(0,0,0) & 28 & 0 & 0 \cr
(1,0,0) & 16 & 6 & 0 \cr
(2,0,0) & 8 & 8 & 1  \cr
(3,0,0) & 4 & 6 & 3  \cr
(0,1,0) & 10 & 6 & 0 \cr
(0,2,0) & 1 & 6 & 1 \cr
(0,3,0) & 1 & 0 & 3 \cr
(1,1,0) & 4 & 7 & 1 \cr
(2,1,0) & 2 & 4 & 3 \cr
(1,2,0) & 1 & 2 & 3 \cr
(0,0,1) & 6 & 5 & 0 \cr
(1,0,1) & 2 & 5 & 1 \cr
(0,1,1) & 0 & 4 & 1}
$$
 \
\ni
Our  values for $b_0$ agree with Joe Harris's Theorem 3.7 of [H] (wherever  ours and his table intersect).

\

\ni
{ \bf \irreducible}
The computation of the multiplicities of the theta-lines in all the above cases
 is contained in the following lemma.
\

\proclaim Lemma \irreduciblelemma . Let $X$ be an irreducible quartic in $V$. Then:
\item\item
a) A theta-line of type $0$ appears in $\t(X)$ with multiplicity 1.
\item\item
b) A theta-line of type $1$ containing a node (resp. a cusp)
appears in $\t(X)$  with multiplicity $2$ (resp. $3$).
\item\item
c) A theta-line containing two nodes  appears in $\t(X)$  with multiplicity $4$.
\item\item
d) A theta-line containing two cusps appears in $\t(X)$  with multiplicity $9$.
\item\item
e) A theta-line containing a node and a cusp 
appears in $\t(X)$  with multiplicity $6$.
\item\item
f) A theta-line of type $1$ containing a tacnode 
appears in $\t(X)$  with multiplicity $4$, unless it is the tacnodal 
tangent, which instead appears with multiplicity $6$.
\item\item
g) A theta-line of type $2$ containing a tacnode 
appears in $\t(X)$  with multiplicity $8$ if it contains a node, and 
$12$ if it contains a cusp.

\pf
Let us suppose first that $X$ has no hyperflexes. Let 
Sing$(X)\subset X$ be its singular locus,  and let $X^{0}=X \backslash {\rm Sing}(X)$.
Consider in $X^{0}\times X^{0}$ the ``tangential correspondence":
$$
Z^0 = \{(x,y): x\ne y \ \hbox{and $y$ belongs to the tangent line to $X$ at $x$}\}
$$
The first projection $Z^0 \to X^{0}$ is a generically $2$ to $1$ map; we let
$Z = \overline{ Z^0} \subset X\times X$ and $\pi:Z \to X$ the first projection.

Now consider a $1$-parameter family of plane quartics $\X \to T$ such that
$X_{0}=X$ and $X_{t}$ is nonsingular for all $t\ne t_0$. We have a corresponding
family of double covers:
$$
\matrix{
\Z&& \cr
\mapdown{\Pi}& \cr
\X& \to & T }
$$
such that $Z_{t} \to X_{t}$ is the tangential correspondence for all 
$t\in T$.
Let $x_0\in X$ be on a theta-line.
In order to compute the multiplicity
of this theta-line in $\t(X)$ we have to compute the number of branch points of
$\Pi_{t}$ which tend to $x_0$ as $t$ tends to $t_0$. 
Let $v$ be a local coordinate in $T$ around $t_0$. We have the following cases:

\

{\it  $x_0\in X^{0}$  is on a theta-line of type} $0$.  
Then $x_0$  is an ordinary branch point of $\pi$ and therefore it is a limit of just
one  branch point of $\Pi_{t}$. This proves (a).

\

{\it  $x_0\in X^{0}$  is on a theta-line of type $1$ containing a node}.

\ni 
  The surface $\Z$ has equation
 $xy-v=0$ locally at $\Pi^{-1}(x_0)$,
 and $\Pi(x,y,v)= (x+y,v)$; the map $\X \to T$ is locally
$(u,v) \mapsto v$ ($u$ being a local coordinate in $X$ around $x_0$).
 Then for a given $v\ne 0$ the branch points of $\Pi_{t}$ are
$(\pm 2\sqrt{v},v)$, so there are two branch points tending to $x_0$ and this proves 
(b) in the nodal case. 

\

{\it  $x_0\in X^{0}$  is on a theta-line of type $1$ containing a cusp}.

\ni Then  
  the surface $\Z$ has equation
 $x^2-y^3+v=0$ locally at $\Pi^{-1}(x_0)$,
and $\Pi(x,y,v)= (y,v)$. For a given $v\ne 0$ the branch points of $\Pi_{t}$ are
$(\zeta_1,v), (\zeta_2,v), (\zeta_3,v)$, where $\zeta_{1},\zeta_{2},\zeta_{3}$ are the  cubic roots of $v$; 
this proves  (b) in the cuspidal case.

\

{\it  $x_0\in X^{0}$  is on a theta-line of type $1$ containing a tacnode and
different from the tacnodal tangent}.

\ni
The surface $\Z$ has equation
 $x^2-y^4+v=0$ locally at $\Pi^{-1}(x_0)$,
and $\Pi(x,y,v)= (y,v)$. 
For a given $v\ne 0$ the branch points of $\Pi_{t}$ are 
$(\eta_i,v)$, $i=1\ldots,4$,  where $\eta_{1},\ldots,\eta_{4}$ are the  
quartic roots of $v$. This proves the first part of (f). 

\

{\it $x_0\in X$  is on a theta-line of type $2$}.

\ni
Assume first that $x_{0}$ is a node and call  the other singular point on 
the theta-line  $p$. Let $\X \to T$ be a family of quartics such 
that $X_{0}=X$ and for each $t\ne t_{0}$ $X_{t}$ has a node $x_{t}$ 
specializing to  $x_{0}$, and no other singularity.  $\X$ has a double 
curve generated by the varying node  $x_{t}$. Let $\nu: \X' \to \X$ be 
the normalization. Then $\X' \to T$ is a family  of generically smooth 
curves. Let $\sigma:\Z' = \X'\times_{T}\Z \to \X'$ and let
$$
\matrix{
\sigma_{0}: & Z'& \to & X' \cr
&&& \mapdown{\nu_{0}} \cr
&&& X}
$$
be the central fiber. Let $\nu_{0}^{-1}(x_{0}) = \{y,y'\}$. Then, 
arguing as before, we deduce that $y$ and $y'$ are each limit of the 
appropriate number $\ell$ (two, three, four) of branch points of $\sigma_{t}$.
Thus $x_{0}$ counts with multiplicity $2\ell$. This takes care of 
theta-lines of type $2$ containing a node. 

The case of a theta-line joining a tacnode and a cusp is treated 
similarly, by considering a degeneration of binodal curves whose 
nodes specialize to the tacnode of $X$. 

There are two remaining cases: a theta-line joining two cusps and the 
tacnodal tangent line. Of course they cannot occur on the same curve.
If $X$ has two cusps and no other singularities and $L$ is the line joining them, then the 
multiplicity $mult_{\t (X)}(L)$ of $L$ in $\t(X)$ is given by the formula
$$
mult_{\t (X)}(L) = 28 - 1-6\cdot 3 = 9
$$
since $X$ has one theta-line of type $0$ and six of type $1$ (see the 
table above).

Similarly, if $X$ has a tacnode and no other singularities and $L$ is 
the tacnodal tangent 
$$
mult_{\t (X)}(L) = 28 - 6 -4\cdot 4 = 6
$$
since such a tacnodal quartic has six theta-lines of type $0$ and four 
of type $1$ (other than the tacnodal tangent).

The careful reader can verify that the same conclusion can be 
achieved by the same analysis applied to other examples.

\

Finally notice that if $X$ has hyperflexes, then it can be obtained 
as a limit of curves without hyperflexes, and having the same type of 
singularities. Since the multiplicities of the components of the 
theta-curves can only increase by specialization and the number of 
theta-lines of type $b_{i}$ is constant, such multiplicities must remain 
unchanged. The conclusion follows.
\qed

\ 

\rk 
The proof shows that theta-lines of the same type, passing through the same type of singularity,
all have the same multiplicity.

 Notice also that most of the argument does not
require
$X$ to be irreducible.

\

\ni
{ \bf \reducible}
Now we describe the theta curve  for those reducible curves in $V$ having finite stabilizer
(which is all we need); we refer to [AF1] for a study of quartics with
infinite stabilizer. In this case, the support of the theta-curve is easy to find, by the
Hurwitz formula, or looking at the dual curves (when
$X$ is a union of conics). The multiplicities of those components of $\t (X)$, which are
not also components of
$X$, are computed just like in \irreduciblelemma (see the above 
remark). 
If $X$ contains a line $L$, then $L\subset \t (X)$; its multiplicity is easily obtained
in each case (see below).

\

\ni
{\it Union of two smooth conics: $X=C_1\cup C_2$}

\ni
1. If $X$ has $4$ nodes, then $X$ is a so-called {\it Split curve} and it is treated in detail in  section
\split  below.

\ni
2. If $X$ has $2$ nodes $n_1$ and $n_2$ and a tacnode $t$, that is, the two conics are tangent in $t$,
then 
$$
\t (X) = L_1 \cup L_2 \cup 4 E \cup 6 L \cup 8 (M_1\cup M_2)
$$
where $L_1$ and $L_2$ are the two lines of type $0$ (corresponding to the two nodes of the dual curve
of $X$, which is also the union of two conics, tangent at one point);
$E$ is the line through the two nodes, $L$ is the tangent line to the conics at $t$
(which is a tacnodal tangent) and $M_i$ joins $t$ and $n_i$ (see \irreduciblelemma ).

\

\ni
{\it Union of a line $L$ and an irreducible cubic $C$  meeting transversally at $n_1,
n_2, n_3$}

\ni
3. If $X$ has only $3$ nodes and no other singularity, then by Hurwitz formula
$$
\t (X) = 2 (\cup _{i=1}^3(M_1^i\cup M_2^i\cup M_3^i\cup M_4^i)) \cup 4L
$$
where $M_j^i$ is a theta-line of type $1$ passing through the node $n_i$ and tangent to $C$.

\ni
4. If $C$ has a node $n$, then
$$
\t (X) = 2 (\cup _{i=1}^3(M_1^i\cup M_2^i)) \cup 4 (\cup _{i=1}^3(E_i)) \cup 4L
$$
where $E_i$ joins $n$ with the node $n_i$.

\ni
5. If $C$ has a cusp $c$, then
$$
\t (X) = 2 (\cup _{i=1}^3(M^i))  \cup 6(\cup _{i=1}^3(E_i))\cup 4L
$$
where $E_i$ joins $c$ with the node $n_i$.

\

\ni
{\it Union of a line $L$ and an irreducible cubic $C$ tangent to $L$ at $t$}

Here we call $n$ the remaining point of intersection of $L$ and $C$; recall that,
for $X$ to be GIT-semistable, it is necessary that $t$ is not an inflectionary point
of $C$ ([GIT] 4.2). Clearly $t$ is a tacnode for $X$.

\ni
6. If $C$ is smooth, then
$$
\t (X) =2 (M_1\cup M_2\cup M_3) \cup 4 (\cup _{i=1}^4(K_i)) \cup 6L
$$
where $n\in M_i$ and $t \in K_j$ (all these lines are of type $1$)

\ni
7. If $C$ has a node $m$
$$
\t (X) =2 M \cup 4E \cup 8F \cup 4 (\cup _{i=1}^2(K_i)) \cup 6L
$$
where $E$ is the line through $n$ and $m$ and $F$ the line through $m$ and $t$.

\ni
8. If $C$ has a cusp $c$ 
$$
\t (X) =  4 K \cup 6E \cup 12F \cup  6L
$$
where $t\in K$ and $E$ is the line through $n$ and $c$ and $F$ the line through $c$ and
$t$.

\

\ni
{\it $X$ contains two distinct lines $L$ and $M$}

Call $n=L\cap M$.

\ni
9. $X= L\cup M \cup C$, with $C$  smooth and transverse to both lines, then
$$
\t (X) = 2 (L_1\cup L_2) \cup 4(E_1\cup E_2\cup E_3\cup E_4\cup L\cup M)
$$
where $ L_i $ contains $n$ and is tangent to $C$, and the remaining are the $6$ lines
joining the $4$ pts of $C\cap (L_1\cup L_2)$.

\ni
10.
$X= L\cup M \cup C$, with $C$  smooth, tangent to $M$ at a point $t$ and meeting $L$ at
$n_1$ and $n_2$.
$$
\t (X) = 2 L_1 \cup 4 L \cup 6M \cup 8(E_1\cup E_2)
$$
where $ L_1 $ contains $n$ and is tangent to $C$, $E_i$ joins $t$ with $n_i$.

\ni
11. 
$X$ is the union of $4$ general lines $L_1,....,L_4$, then 
$$
\t(X) = 4(L_1\cup....\cup L_7)
$$
in fact, of the $6$ nodes of $X$ there are $3$ pairs not lying on any line of $X$;
$L_5, \ L_6,\ L_7$ are the lines joining such pairs.

\ni
{\bf \tacnodal}
The behaviour 
of theta-curves
of quartics in $V$ with respect to the action of $G=PGL(3)$ on $P_{28}$  is described by the following

\proclaim Lemma \key . Let $X\in V$ be such that $Stab _G(X)$ is finite; let $\T =\t (X)$.
Then either $\T$ is GIT-stable (as a point in $P_{28}$, acted on by $G$) or
$X$ is a tacnodal curve, more precisely, $X$ is one of the following types

{\it 1. irreducible with one tacnode and no other singularity;

 2. irreducible with one tacnode and one node;

3. irreducible with one tacnode and one  cusp;

4.  the union of two smooth conics tangent at one point and intersecting transversally
in 2 other points;

5. the union of a line $L$ and a smooth cubic $C$, with $C$ tangent to $L$ at a non-inflectionary
point;

6. the union of a line $L$ and an irreducible, nodal cubic $C$, with $C$ tangent to $L$ at a
non-inflectionary point;

7. the union of a line $L$ and an irreducible, cuspidal cubic $C$, with $C$ tangent to $L$ at a
non-inflectionary point.

8. the union of a smooth conic $C$ and two lines $L$ and $M$ such that $M$ is tangent
to $C$ and $L$ is transverse.}

\

\pf
By \mumford all semistable points in $P_{28}$ are actually stable, and, if $\T$ is not stable, then it has
either a line of multiplicity at least $10$ or a point of multiplicity at least $19$.
We refer  to [AF1] for the list of quartics with infinite stabilizer.
Looking at the description in \irreducible and \reducible, we see that the only quartics
in $V$ with finite stabilizer, having  a theta-line of multiplicity  at least $10$ are those with
a tacnode and a cusp, i.e. cases  3 and 7 above. Analogously, the only theta-curves having
a point of multiplicity at least $19$ are the theta-curves of tacnodal curves, where the
tacnode of $X$ becomes a point of multiplicity $22$ for $\t (X)$.
\qed

\

\ni
\rk From the results of this section, it follows easily that if two quartics in $V$ have the same theta-curve,
then they have the same  singularities and decomposition in irreducible components of the same type.
Here are some relevant examples, for all of which we restrict our attention to curves in $V$.

\

(a) $X$ is smooth if and only if $\t (X)$ is made of 28 distinct lines.

\

Notice now that a line of multiplicity 2 is a theta-line of type 1 passing through a node, thus

(b) $X$ has one node and no other singularity if and only if $\t (X)$ contains 6 incident lines of
multiplicity 2 and no line of higher multiplicity. The node is determined by the incidence point of such 6
lines.

\

If $\t(X)$ contains lines of multiplicity 1, then $X$ does not contain any line as a component
(a cubic does not have any bitangent line!). In such a case, the multiplicity of a line of $\t (X)$ uniquely
determines the type of singularities that it contains.
We can apply this to the following examples

(c) $X$ has 2 nodes and no other singularity if and only if $\t (X)$ contains some line of multiplicity 1,
one line of multiplicity 4 and two sets of 4 incident lines of multiplicity 2. The two nodes are
determined by the 2 incidence points of the two sets above.

(d) $X$ has 3 nodes (and no other singularity) if and only if $\t (X)$ contains some line of multiplicity 1
and 3
non-concurrent lines of multiplicity 4. The nodes are the 3 verteces of the triangle formed by such 3 lines.

(e) $X$ is the union of 2 smooth conics meeting transversally ($X$ is a ``split" curve, see below)
if and only if $X$ has 4 non-colinear nodes and no other singularity, if and only if $\t(X)$
contains some line of multiplicity 1
and a set of 6 lines of multiplicity 4 containing 4 triples of incident lines.
The nodes are the intersection points of such 4 triples.

\

\

\

\cl{4. SPECIAL CASES}

\

\ni
{\bf\split}
We prove here the theta-property  \ for two types of singular quartics.

\

\df Let $X$ be a plane quartic having two  irreducible components: $X= C_1\cup C_2$
such that $C_1$ and $C_2$ are smooth conics  meeting transversally. Then $X$ is called a {\it
split curve}.

\

Thus a split curve $X=C_1\cup C_2$  has $4$ nodes 
and no other singularities, in particular it is
GIT-stable. Denote by $\{ n_1,n_2,n_3,n_4\} = C_1\cap C_2$ the nodes of $X$.
Let $E_{i,j}$ be the line through $n_i$ and $n_j$, it is a theta-line of type $2$ for $X$.
The theta-curve $\T$ of $X$ is the following
$$
\T = L_1\cup L_2\cup L_3\cup L_4\cup 4\bigl(\cup_{i,j}E_{i,j}\bigr)
$$
where each of the 6 lines $E_{i,j}$ counts with multiplicity $4$ by \irreduciblelemma ;
the remaining 4 lines $L_i$ are of type 0, and they are determined by the intersection
of the dual curves $C_1^*$ and $C_2^*$ in the dual plane.

We shall now prove that a split curve is determined, among all curves in $V$, by its theta-curve.

\

\proclaim Proposition \splitpr .
Split curves have the theta-property.

\pf
We keep using the  notation above. Let $Y\in V$ be a curve such that $\t (Y) = \T$, where $\T$ is
the theta-curve of a split-curve $X$.
By case (e) of the remark at the end of the previous section,  $Y$ is also a split curve and has the same
nodes of $X$.

 Let $P\cong \P ^1$ be
the pencil of all conics passing through 
$n_1,....,n_4$.
Of course  $Y$ and $X$ are   union of two conics   belonging to $P$.
Now consider the line $L_1$; in $P$ there are exactly two conics that are tangent to $L_1$, in
fact $P$ cuts on $L_1\cong \P ^1$ a linear series of degree $2$ and dimension $1$ which
has $2$ ramification points, corresponding to those conics of $P$ that intersect $L_1$ in a unique
point. We conclude that such conics are the original $C_1$ and $C_2$ and that $X=Y$.
\qed

\

\ni
{\bf \trinodal}
Now we treat a second special case: that of irreducible quartics with three nodes.
We call such curves {\it trinodal}. If $X$ is a trinodal curve, denote by $n_1,n_2,n_3$ its nodes;
let $\T$ be its theta-curve, then $\T$ contains the three lines $\overline{n_i,n_j}$ with
multiplicity $4$, three pairs of theta-lines of type 1, $L_i,M_i$ with $i=1,2,3$ such that $L_i\cap
M_i = n_i$ (hence each of these six lines appears with multiplicity $2$ in $\T$) and four
theta-lines of type 0.

\proclaim Proposition \trinodalpr . A trinodal curve has the theta-property.

\pf
Let $X$ and $\T$ be as above. Let $Y\in V$ be such that $\t (Y)=\T$; then,
by the remark at the end of \tacnodal   $Y$ is also a trinodal
curve with the same nodes of $X$.
Let $\g :\Pl \la \Pl $ be the quadratic transformation with base points $n_1,n_2,n_3$
(that is, $\g$ is the birational map associated to the linear system  of all conics
passing through $n_1,n_2,n_3$).
Then $\g$ maps the set of all quartics that are singular at $n_1,n_2,n_3$ bijectively to the set
of all conics in the image $\Pl$.
Thus $X$ and $Y$ are mapped to smooth conics $C$ and $D$, and $X$ and $Y$ are distinct if and only
if
$C$ and $D$ are distinct. Consider now the six theta-lines of type 1, $L_i,M_i$ for $i=1,2,3$;
their images under $\g$ are six lines that are tangent to $C$ and $D$. But now given $6$  
lines in $\Pl$, there exists at most one smooth conic that is tangent to all of them
(this is seen easily by passing to the dual plane $\Pld$, where the $6$ lines correspond to
$6$ points and the dual of a smooth conic tangent to those $6$ lines is a smooth conic passing 
through those $6$ points).
We hence conclude that $C=D$ and therefore $X=Y$. \qed

\

\

\cl{5. THE MAIN RESULT}

\

\ni
{\bf \binodal } Our main goal is to prove Theorem \main .
The two next statements, \binodalpr \  and \nodalpr  are completely independent but they
are    proved in a
very similar fashion. We include them for the sake of completeness, and to follow the pattern
started in the previous section, proving the theta property for nodal quartics with decreasing number of nodes.

We shall say that a plane curve is binodal if it has two ordinary double 
points and no other singularity.

\proclaim Proposition \binodalpr . A general binodal quartic has the theta-property.

\pf
By contradiction, assume that the statement is false. Then we can construct a family
of curves $\X \la T$ over a smooth curve $T$ with a marked point $t_0$, with the following properties: there
are two fixed points
$n_1$ and
$n_2$ in $\Pl$ such that
if $t\neq t_0$, the fiber $X_t$ has two nodes at $n_1$ and $n_2$ and no other singularities; the special
fiber $X_0$ is trinodal so that necessarily two of the nodes of $X_0$ are in $n_1$ and $n_2$. 
We also assume that each theta-line of $X_0$ meets it in two distinct points (this condition is satisfied for a general
trinodal curve). Finally, we require that $X_t$ does not have the 
theta-property, so that if we denote by $\T _t =\t (X_t)$, there exists a second family
$\Y \la T$ with $Y_t\in V$ if $t\neq t_0$, $X_t\neq Y_t$, and such that $\t (Y_t) = \T _t$ for every $t$.

Notice that $Y_t$ is necessarily a binodal curve with nodes in $n_1$ and $n_2$, by the analysis
in section 3. Consider now $Y_0$; there are of course two possibilities: (a) $Y_0 \in V$;
(b) $Y_0 \not\in V$.

Case (a).

By construction $\t (Y_0) = \t (X_0)$; this implies that $Y_0=X_0$ because
$X_0$ has the theta-property, by \trinodalpr .
Thus the two families $\X$ and $\Y$ have the same special fiber. We shall now prove that this
forces $X_t = Y_t$ for every $t$.
Denote by $W$ the $8$-dimensional linear subspace of $\P ^{14}$ consisting of all quartics that
have a singular point at $n_1$ and $n_2$. Our $\X$ and $\Y$ are thus families of curves in $W\cap
V$. Let $\t _W:W\cap V \la P_{28}$ be the restriction of $\t$, we shall prove the following claim:

\ni
{\it $\t _W$ is an immersion at $X_0$}.

It is clear that this implies that $X_t = Y_t$, which is a contradiction.

To prove the claim, let us first describe,
\`a la Zariski ([Z]),  the tangent space to $W$ at any point $X$. 
 This is identified to the linear
series cut on $X$ by curves of degree $4$, doubled at the nodes of $X$, that is,
let $\nu : X^{\nu }\la X$ be the normalization of $X$, then

$$
T_XW \cong H^0(X^{\nu} , \nu ^*(\O _X(4)\otimes \O_X(-2n_1-2n_2)))
$$
let $\L :=\nu ^*(\O _X(4)\otimes \O_X(-2n_1-2n_2))$, so that $\deg \L = 8$.
Now, let $L_1,....,L_4$ be the $4$ lines of type $0$ contained in $\t (X_0)$,
and let $M$ be any theta-line of type $1$;
 denote by $W_{L_1,....,L_4}^M$ the variety of curves in $W$ that are bitangent to $L_1,....L_4$,
and tangent to $M$.
Then we have (again using Zariski's theory)
$$
T_{X_0}W_{L_1,....,L_4}^M\cong H^0(X_0^{\nu}, \L\otimes \nu ^*\O(-p_1-p_2-....-p_8-q))
$$
where $p_1,....p_8$ are the $8$ points of intersection of
$X_0$ with $L_1\cup....\cup L_4$ (with $p_i\neq p_j$ by our assumption that $X_0$ has no hyperflexes), and $q$
is the smooth point of intersection of
$M$ and
$X$. Therefore, since $\deg \L\otimes \O(-p_1-p_2-....-p_8-q)<0$, we obtain that
$T_{X_0}W_{L_1,....,L_4}^M=0$.

Now notice that 
$$
T_{X_0}\t _W^{-1}(\t(X_0)) \subset T_{X_0}W_{L_1,....,L_4}^M
$$
hence a fortiori $T_{X_0}\t _W^{-1}(\t(X_0)) =0$, so that $\t _W$ is an immersion at $X_0$.

This concludes the proof of our proposition in case (a).

Case (b).

Suppose now that $Y_0\not\in V$. Thus $Y_0$ is either GIT-unstable, or it is a double
conic (and hence it has infinite stabilizer).
We now apply the GIT semistable replacement property (see \GITSRP ); this gives, up to
replacing
$T$ with a finite covering ramified only over $t_0$, a new family $\Z \la T$ of plane
quartics,  and a morphism $\g : T-\{ t_0\} \la G$
such
that away from $t_0$ we have
$$
Z_t=Y_t^{\g (t)},
$$
and we can assume that $Z_0$ is GIT-semistable with finite stabilizer.

We thus have that $\t(Z_t)=\T_t ^{\g (t) }$ away from $t_0$. 

Recall now that by \key  we have that if $\T \in P_{28}$ is the theta-curve of a reduced quartic having
at most ordinary nodes or ordinary cusps as singularities, then $\T$ is GIT-stable.
And the curves in $V$ having GIT-unstable theta curve are exactly the tacnodal curves.

 Consider now the 
(rational) quotient map
$$
q:P_{28}\la P_{28}//G
$$
 regular
on the set of theta-curves of smooth, nodal and cuspidal quartics.

We have, by construction, $q(\T _t) = q(\t (Z_t))$ for every $t\neq t_0$.
Hence $\t (Z_0) $ lies in the closure of the orbit of $\T _0$.
Denoting  $\T _0'=\t (Z_0)$  we write:
$$
\T _0'\in \overline{O_G(\T_0)}.
$$
There are, of course, two possibilities: $\T_0'$ may or may not be in $O_G(\T_0)$.

If $\T '_0 = \T _0^g$ for some $g\in G$, then $Z_0$ is a trinodal curve (by the analysis in Section 3) and,
as such, it has the theta-property (see \trinodalpr ). In particular, since $\t (X_0^g) = \T
'_0$ we obtain that $Z_0=X_0^g$. But then, when applying the GIT-semistable replacement
property to $\Y \la T$, we can choose $X_0$, rather than $Z_0$, as a semistable
replacement for $Y_0$ (\GITSRP ). Thus, after acting  again with $G$, we are back to case (a).

Finally, suppose that $\T _0'$ is not in the $G$-orbit of  $\T _0$. Since
$\T_0$ is GIT-stable (by \key ),
$\overline{O_G(\T_0)}-  O_G(\T_0)$ only contains GIT-unstable points.
Therefore $Z_0$ must be a tacnodal curve with finite stablizer.
$Z_0$ must then be one of   the eight types listed in \key  .

Let us start with case 1 (of Lemma \key ) and show that it cannot occur. Recall that the
theta-curve
$\T _0'$ of such a tacnodal quartic contains $6$ distinct lines of multiplicity $1$
(the $6$ theta-lines of type 0). We are saying that $\T
_0'$ is in the orbit closure of the theta curve $\T _0$ of a trinodal quartic. We know that $\T
_0$ contains only
$4$ lines of multiplicity $1$ (the $4$ theta-lines of type $0$), and hence any configuration of
lines in its orbit closure has at most $4$ lines of multiplicity $1$. Thus case 1 does not occur.

Instead of using a similar (even though slightly more complicated) argument 
to show that cases
2, 3, 4 cannot occur, 
we observe that all such curves, and curves of type 1 
as well, contain in their orbit closure 
reducible quartics that are union of two
 bitangent conics (that is, quartics $C_1\cup C_2$ with $C_1$  and $C_2$ smooth conics
having contact of order $2$ at exactly two points;
 see [AF2]). Such reducible quartics are 
GIT-semistable; moreover   for each curve of type 2, 3, 
or 4 there is a curve of type 1 with the same pair of bitangent conics 
in the orbit
closure. In other words, every curve of type 2, 3, or 4 is 
identified with one of type 1 in the GIT-quotient of $\Bbb P^{14}$. 
In particular,
when we applied the GIT-semistable replacement property to obtain
$Z_0$ we were free to assume that $Z_0$ were a tacnodal curve of type 1, rather than
one of type 2, 3 or 4. Since we just showed that 
$Z_0$ cannot be of type 1, we also get that $Z_0$ cannot be of type 2, 3 or 4.

Consider now case 5: $Z_0 = L\cup C$.
Then there are no lines of multiplicity $1$, since $Z_0$ has no theta-lines of type $0$;
$\T _0'$ contains 
$L$ with multiplicity $6$, one set of $4$ lines
of multiplicity
$4$ all meeting $L$ at the same point $t$ (the tacnode,  where $C$ is tangent to $L$) and $3$
lines of multiplicity $2$ all meeting $L$ at the same point $n$  (the remaining point of intersection of
$C$ and
$L$). Therefore $t$ is a point of multiplicity $22$ and $n$ of multiplicity $12$; every other
point of $\T _0'$ has multiplicity less than $12$. But now consider $\T _0$, which contains $3$
points of multiplicity $12$ (corresponding to the $3$ nodes of $X_0$). Then a curve in the orbit
closure of
$\T _0$ must contain either at least $3$ points of multiplicity at least $12$, or at least one
point of multiplicity at least
$24$. Neither condition is satisfied by the theta-curve of our $Z_0$, hence case 5 does not occur.
(Notice that the same argument would work to show that cases 1,2,3 are not possible).

We treat the three remaining cases by noticing that they are both identified in the quotient $\P
^{14} // G$ to curves of type 5, since all such curves contain in their orbit closure
the GIT-semistable quartics given by a smooth conic union 2 tangent lines (see [AF2]).
Hence, arguing as above, when applying the GIT-semistable replacement property to obtain $Z_0$ we
could assume that $Z_0$ was a  curve of type 5, rather than one of type 6, 7 or 8.

\qed

\proclaim Proposition \nodalpr . A general nodal quartic has the theta-property.

The argument is completely similar, we just need to replace $W$ with the linear space of quartics
having a node in a fixed point $N_1$. This amounts to making a few obvious changes to the proof of
Case (a) (consisting in replacing  the line $M$ with all theta-lines of type 1 of $X_0$). The
rest remains the same.

\

\ni
{\bf\smooth }
We now prove that a general quartic is uniquely recovered from its 28 bitangents.
\proclaim Theorem \main . A general smooth plane quartic has the theta-property.

\pf
We proceed like in the proof of \binodalpr .
By contradiction, assume that for $X$ varying in an open subset of $\Vs$ there exists $Y\in \Vs$
such that $\t(X)=\t(Y)$ and $X\neq Y$.

Then we can construct a family of  quartics $\X\la T$ over a smooth curve $T$
with a marked point $t_0\in T$ satisfying the following assumptions: we fix two distinct lines
$E$ and $F$ and we assume that, for every $t\neq t_0$, the fiber $X_t$ is smooth
and has $E$ and $F$ as bitangents; the fiber over $t_0$,
$X_0$, is a trinodal curve, having all nodes on $E$ and $F$, with one node at the point of
intersection of $E$ and $F$
(so that $E$ and $F$ are theta-lines of type $2$ for $X_0$). 
Just like in the proof of \binodal  it is convenient to
assume that each theta-line of $X_0$ meets it in two distinct points
(this condition is satisfied for a general trinodal curve).
Such degenerations do exist: in fact denote by $V_{E,F}$
the subvariety of $V$ parametrizing curves having $E$ and
$F$ as theta-lines. Then 
$V_{E,F}$ is irreducible and its general
element is a smooth quartic. To see that, consider the dominant, rational map
$$
c:V_{E,F} \la Sym^2E\times Sym ^2F
$$
mapping a curve $X$ to its two pairs of points of intersection with $E$ and $F$.
The general fiber of $c$ is irreducible, being just the intersection of $V$ with a linear subspace of
dimension 6 of $\P^{14}$. Hence $V_{E,F}$ is irreducible of dimension $6+4=10$. The fact that the general
element of $V_{E,F}$ is a smooth curve can be seen by a standard tangent space argument (see below).

Finally, we assume that there exists another family
$\Y
\la T$ such that for every
$t\neq t_0$
$Y_t$ is in $V$, $X_t\neq Y_t$  and $\t (X_t)=\t (Y_t)$.

Notice that if $t\neq t_0$, $Y_t$ must itself be smooth, since its theta-curve is made of 28
distinct lines.

Our $X_t$, $X_0$ and $Y_t$ are in $V_{E,F}$.
What about $Y_0$?  There are two possibilities:

(a) $Y_0 \in V$ and hence $Y_0 \in V_{E,F}$

(b) $Y_0 \not\in V$

Case (a).

Then $X_0=Y_0$, because $X_0$ has the  theta-property, by \trinodal  Just as in \binodalpr , we will be
done after showing that the restriction 
$$
\theta _{E,F}: V_{E,F}\la P_{28}
$$ 
of $\theta $ to $V_{E,F}$,  is an immersion at $X_0$.

Notice that, if $X\in V_{E,F}$  is smooth and $X$ intersects $E$ and $F$ in $p_1, ....,p_4$,
then
$$
T_XV_{E,F}\cong H^0(X,\O _X(4)\otimes \O _X(-p_1- ....-p_4))
$$
which has dimension $10$ by Riemann-Roch.
If, instead, we consider our $X_0$, then we have that
$$
T_{X_0}V_{E,F}\subset H^0(X_0^{\nu },\nu ^*\O _{X_0}(4)\otimes \O _{X_0}(-n_1-n_2-n_3))
$$
where $\nu:X_0^{\nu }\la X_0$ is the normalization.
Denote by $\L =\nu ^*\O _{X_0}(4)\otimes \O _{X_0}(-n_1-n_2-n_3)$.

Now, the theta-curve of $X_0$ contains $4$ lines $L_1,....,L_4$ of type $0$,
meeting $X_0$ in $8$ smooth points $p_1,....,p_8$ ($2$ distinct on each line, since we picked $X_0$ without
hyperflexes) and
$6$ lines 
$M_1,....M_6$ of
type
$1$ meeting $X_0$ in $6$ smooth points $q_1,....,q_6$ ($1$ on each line).
Denote by
$V_{E,F,L_1,L_2,L_3,L_4}^{M_1,....M_6}$ the locus of curves in $V_{E,F}$ that are bitangent to
all $L_i$'s and tangent to all $M_j$'s, we have
$$
T_{X_0}V_{E,F,L_1,L_2,L_3,L_4}^{M_1,....M_6}\subset H^0(X_0^{\nu},\L\otimes \nu ^*\O_{X_0}(-\sum _1^8p_i-\sum
_1^6q_j))
$$
and the right hand side vanishes, since the line bundle has negative degree.
On the other hand it is clear that
$$
T_{X_0}\t _{E,F}^{-1}(\t (X_0) \subset T_{X_0}V_{E,F,L_1,L_2,L_3,L_4}^{M_1,....M_6}=0
$$
hence $\t _{E,F}$ is an immersion at $X_0$ and we are done.

Case (b).
Here we just repeat, word by word,  the argument for Case (b) in \binodalpr !
\qed

\

\

\ni
{\bf Aknowledgements.} We benefitted from  conversations with Igor Dolgachev, Lawrence Ein, Joe
Harris and Sandro Verra. We are  very grateful to Paolo Aluffi for promptly and very clearly answering, 
by e-mail, all of our many  questions.

\vfill\eject
\null

\cl{REFERENCES}

\

\ni
[AF1] Aluffi P. - Faber C.: {\it Plane 
curves with small linear orbits II.}  International J. of Math, to appear.

\

\ni
[AF2] Aluffi P. - Faber C.: {\it Linear orbits of arbitrary plane 
curves.}  Mich. J. of Math (W.Fulton issue), to appear.

\

\ni
[ACGH] Arbarello - Cornalba - Griffiths - Harris: {\it Geometry of Algebraic Curves, I} Grund. series,
vol. 267, Springer.

\

\ni
[E-Ch] Enriques F. - Chisini O.: {\it Teoria geometrica delle equazioni e delle
 funzioni algebriche}, vol. I. Zanichelli, Bologna 1929.
 
 \

\ni
[H] Harris J.: {\it Theta-characteristics on algebraic curves.} Trans. AMS 271 N.2 June 1983

\

\ni
[K-W] Krazer A. - Wirtinger W.: {\it Abelsche Funktionen und Allgemeine 
 Thetafunktionen}. Enc. d. Math. Wiss. II B 7, Leipzig 1921.

\

\ni
[GIT] Mumford D. - Forgarty J. - Kirwan F.  {\it Geometric Invariant Theory.} (Third edition)
E.M.G. 34 Springer 1994.

\

\ni
[S] Segre C.: {\it Le molteplicit\`a nelle intersezioni delle curve piane algebriche con 
alcune applicazioni ai principii della teoria di tali curve.}
 Giornale di Matematiche 36 (1898), 1-50. Opere, vol. I, 380-429.
 
 \
 
\ni
[W] Walker R.J.: {\it Algebraic Curves.} Princeton U.P. 1950.

\

\ni
[Z] Zariski O. {\it Dimension-theoretic characterization of maximal irreducible  algebraic system of plane nodal curves curves of a
given order $n$ and with a given number $d$ of nodes.} \ Am. J. Math. 104 (1982), pp. 209-226.

\

\ni
Address of Authors:

\ni
L.Caporaso: Universit\`a degli Studi del Sannio, Benevento, Italy

and Massachusetts Institute of Technology, Cambridge, MA, USA

caporaso@math.mit.edu

\

\ni
E.Sernesi: Universit\`a  Roma Tre, Roma, Italy

sernesi@matrm3.mat.uniroma3.it

\


\end